\documentclass[preprint,12pt]{elsarticle}

\usepackage{amssymb}
\usepackage{geometry}
\usepackage{epstopdf} 

\journal{ Journal of Applied Statistics}

\begin{document}

\begin{frontmatter}


\title{On the restricted almost unbiased Liu estimator in the Logistic regression model}


\author{Jibo Wu$^{*1}$, Yasin Asar$^{2}$ and M. Arashi$^{3}$}

\address{ $^1$Key Laboratory of Group  \& Graph Theories and Applications, Chongqing University
of  Arts and Sciences,  Chongqing, 402160, China\\
${^2}$Department of Statistics,
Necmettin Erbakan University,
Konya, TURKEY\\
${^3}$ Department of Statistics, School of Mathematical Sciences, Shahrood University of Technology, Shahrood, Iran
}
\cortext[*]{Corresponding author (Email: linfen52@126.com)}

\begin{abstract}
It is known that when the multicollinearity exists in the logistic
regression model, variance of maximum likelihood estimator is
unstable. As a remedy, in the context of biased shrinkage ridge
estimation, Chang (2015) introduced an almost unbiased Liu estimator
in the logistic regression model. Making use of his approach, when some
prior knowledge in the form of linear restrictions are also available, we introduce a restricted almost unbiased Liu estimator in the logistic regression model.
Statistical properties of this newly defined estimator are
derived and some comparison result are also provided in the form of theorems. A Monte
Carlo simulation study along with a real data example are given
to investigate the performance of this estimator.
\end{abstract}

\begin{keyword}
Almost unbiased Liu estimator; Eigenvalue; Liu estimator; Mean squared error matrix; Restricted almost unbiased Liu estimator\\\vspace{0.3cm}
\textit{AMS Subject Classification}: Primary 62J02, Secondary 62J07
\end{keyword}


\end{frontmatter}


\section{Introduction }
\noindent Consider the following logistic regression model
\begin{equation}
y_i=\pi_i+\varepsilon_i,\ i=1,\ldots,n,
\end{equation}
where
\begin{equation}
\pi_i=\mbox{Pr}(y_i=1)=\frac{exp(x_i'\beta)}{1+exp(x_i'\beta)}, \ i=1,2,\ldots,n
\end{equation}
is the expectation of $y_i$ when the ith value of the dependent variable is Bernoulli with parameter $\pi_i$,
$\beta=(\beta_0,\beta_1,\ldots,\beta_p)'$ denotes the unknown $(p+1)$-vector of regression coefficients, $x_i=(1,x_{1i},x_{2i},\ldots,x_{pi})'$ denotes the ith row of $X$, the $n\times(p+1)$ data matrix, and $\varepsilon_i$s are independent random errors, each having zero mean and variance $w_i=\pi_i(1-\pi_i)$.

In the estimation process of coefficient $\beta$, one often uses the maximum likelihood (ML) approach. Making use of the iteratively weighted least squares algorithm, the MLE is obtained as
\begin{equation}
\hat{\beta}_{MLE}=(X'\hat{W}X)^{-1}X'\hat{W}Z,
\end{equation}
where $Z=(Z_1,\ldots,Z_n)'$, with $Z_i=\log(\hat{\pi}_i)+\frac{y_i-\hat{\pi}_i}{\hat{\pi}_i(1-\hat{\pi}_i)}$ and $\hat W=\mbox{Diag}(\hat{\pi}_i(1-\hat{\pi}_i))$.

In the context of linear regression model, when the multicollinearity exists, the ordinary least squares (OLS) estimator is no longer an efficient estimator. To overcome the problem of multicollinearity, new biased shrinkage estimators have been proposed. As such, Liu (1993) proposed the Liu estimator, the almost unbiased Liu estimator proposed by Akdeniz and Ka\c{c}{\i}ranlar (1995), the restricted Liu estimator proposed by Ka\c{c}{\i}ranlar et al. (1999), the restricted almost unbiased Liu estimator introduced by Xu and Yang (2011a,b) and more recently Asar et al. (2015) developed a two-parameter restricted Liu estimator.

In the logistic regression model, when the explanatory variables are highly correlated, the variance of the MLE becomes inflated. Dealing with such a  problem, some estimators are proposed. One of them is the ridge estimator of Schaefer et al. (1984). When some prior information, in the form of restrictions on regression coefficients, are available, Duffy and Santner (1989) suggested to use the restricted MLE (RMLE).

Now, suppose we are provided with some prior information about $\beta$ in the form of linear restrictions
\begin{equation}\label{restriction}
H\beta=h,
\end{equation}
where $H$ is a $q\times (p+1)$ ($q\leq p+1$) known matrix and $h$ denotes a $q\times 1$ vector of pre-specified known constants.

Considering such a case, Duffy and Santner (1989) defined the RMLE given by
\begin{equation}
\hat{\beta}_{RMLE}=\hat{\beta}_{MLE}-C^{-1}H'(HC^{-1}H')^{-1}(H\hat{\beta}_{MLE}-h),
\end{equation}
where $C=X'\hat{W}X$.

Mansson et al. (2011) defined the Liu estimator (LE) of $\beta$ in the logistic regression model given by
\begin{equation}
\hat{\beta}_{LE}=(C+I)^{-1}(C+dI)\hat{\beta}_{MLE},
\end{equation}
where $0<d<1$ is the biasing parameter.

\c{S}iray et al. (2015) combined the LE and RMLE, to introduce a restricted Liu estimator (RLE) in the logistic regression model as follows
\begin{equation}
\hat{\beta}_{RLE}=(C+I)^{-1}(C+dI)\hat{\beta}_{RMLE},
\end{equation}
where $0<d<1$ is the biasing parameter.

In order to reduce the bias of LE, Chang (2015) proposed an almost
unbiased Liu estimator (AULE) which has form
\begin{equation}
\hat{\beta}_{AULE}=[I-(1-d)^{2}(C+I)^{-2}]\hat{\beta}_{MLE}
\end{equation}
In this paper, using the latter result, a restricted almost unbiased Liu estimator will be exhibited in the logistic regression model when both multicollinearity and restrictions are present.

 We organize the paper as follows: Section 2 contains the introduction of the new estimator along with some preliminary lemmas.
 Comparisons between this estimator with other existing ones are considered in section 3, while a numerical example is conducted in section 4. The work is concluded in section 5.

\section{Proposed Estimator}
In the same fashion as in \c{S}iray et al. (2015), the AULE and RLE
will be combined to obtain a new estimator namely restricted almost
unbiased Liu estimator (RAULE) with form
\begin{equation}
\hat{\beta}_{RAULE}=[I-(1-d)^{2}(C+I)^{-2}]\hat{\beta}_{RMLE},
\end{equation}
where $C=X'\hat{W}X$ and $0<d<1$ is the biasing parameter.\\ Letting $L_d=[I-(1-d)^{2}(C+I)^{-2}]$, the RAULE can be expressed as
\begin{equation}
\hat{\beta}_{RAULE}=L_d\hat{\beta}_{RMLE}.
\end{equation}
For comparison sake and in order to derive characteristics of the RAULE, we need the following lemmas.

\noindent \textbf{Lemma 2.1.} (Shi, 2001) Under the assumptions of section 1, the following matrix is nonnegative definite and has rank $p+1-q$.
$$
A=C^{-1}-C^{-1}H'(HC^{-1}H')^{-1}HC^{-1}
$$
\noindent \textbf{Lemma 2.2.} (Baksalary and Kala, 1983) Suppose
that $M$ be a nonnegative definite matrix and $\alpha$ be a vector,
then
$$
M-\alpha\alpha'\geq 0\Leftrightarrow \alpha'M^+\alpha\leq 1,\ \alpha
\in\Re(M)
$$
where $M^+$ denotes the Moore-Penrose inverse of $M$.

\noindent \textbf{Lemma 2.3.}(Wang, 1994) Suppose that $M$ be a
positive definite matrix and $N$ be a nonnegative definite matrix,
then
$$
M-N\geq 0\Leftrightarrow\lambda_{max}(NM^{-1})\leq 1
$$

\noindent \textbf{Lemma 2.4.}(Wang, 1994) Suppose that both $M$ and $N$ are nonnegative definite matrices, then
$$
M-N\geq 0\Leftrightarrow\lambda_{max}(NM^{-})\leq 1,\
\Re(N)\subset\Re(M)
$$
where $\lambda_{max}(NM^{-})\leq 1$ is invariant of the choice of
$M^-$, and $M^-$ stands for the generalized inverse of $M$.

In the forthcoming section, we will be deriving some properties of the proposed estimator and compare it with some existing competitors.
\section{Properties \& Comparison}
\noindent
Let $\hat\theta$ be an estimator of the parameter $\theta$. The matrix mean squared error (MMSE) of $\hat{\theta}$ is defined by
\begin{eqnarray}
MMSE(\hat{\theta})=E(\hat{\theta}-\theta)(\hat{\theta}-\theta)'=Cov(\hat{\theta})+Bias(\hat{\theta})Bias(\hat{\theta})',
\end{eqnarray}
where $Cov(\cdot)$ is the covariance matrix and
$$
Bias(\hat{\theta})=E(\hat{\theta})-\theta.
$$
The scalar mean squared error (MSE) of $\hat{\theta}$
is defined as
 \begin{eqnarray}
 MSE(\hat{\theta})=\mbox{tr}[MMSE(\hat{\theta})]=\mbox{tr}(Cov(\hat{\theta}))+Bias(\hat{\theta})'Bias(\hat{\theta}).
 \end{eqnarray}
For two given estimators $\hat{\theta}_{1}$ and $\hat{\theta}_{2}$ of $\theta$,
the estimator $\hat{\theta}_{2}$ is said to be superior to estimator
$\hat{\theta}_{1}$ in the sense of MMSE criterion, if and only if
\begin{equation}
\Delta
(\hat{\theta}_{1},\hat{\theta}_{2})=MMSE(\hat{\theta}_{1})-MMSE(\hat{\theta}_{2})\geq
0
\end{equation}

Now, we derive the MMSE of the new estimator (RAULE).

\noindent \textbf{Proposition 3.1.} Under the assumptions of the logistic regression model (1), when the restrictions (\ref{restriction}) hold, the MMSE of the new estimator is given by
\begin{eqnarray}
MMSE(\hat{\beta}_{RAULE})&=&L_dAL_d+\beta'(L_d-I)'(L_d-I)\beta,
\end{eqnarray}
where  $A=C^{-1}-C^{-1}H'(HC^{-1}H')^{-1}HC^{-1}$.\\
\noindent Proof: Using Eq. (9), the covariance and bias of the new estimator are respectively evaluated as
\begin{equation}
Cov(\hat{\beta}_{RAULE})=L_dAL_d
\end{equation}
\begin{equation}
Bias(\hat{\beta}_{RAULE})=(L_d-I)\beta
\end{equation}
Then we obtain
\begin{eqnarray}
MMSE(\hat{\beta}_{RAULE})&=&L_dAL_d+\beta'(L_d-I)'(L_d-I)\beta
\end{eqnarray}
\noindent \textbf{Lemma 3.2.} Under the assumptions of the logistic regression model (1), MMSE of the MLE,
AULE and RMLE are respectively given by
\begin{eqnarray}
MMSE(\hat{\beta}_{MLE}) =C^{-1}
 \end{eqnarray}
\begin{eqnarray}
MMSE(\hat{\beta}_{AULE})
=L_dC^{-1}L_d+\beta'(L_d-I)'(L_d-I)\beta
 \end{eqnarray}
 \begin{eqnarray}
MMSE(\hat{\beta}_{RMLE}) =ACA
 \end{eqnarray}

In the following result, we will be presenting the necessary and sufficient
conditions for the new estimator to be superior to the RMLE in the MMSE
sense.

\noindent\textbf{Theorem 3.3.} Under the assumptions of the logistic regression model (1), assume
$\lambda_{max}(L_dAL_d(ACA)^-)\leq 1$ and $\Re(L_dAL_d)\in\Re(ACA)$.
Then, the RAULE is superior to the
RMLE in the MMSE sense if and only if
\begin{eqnarray}
b_1'(ACA-L_dAL_d)^+b_1\leq1,\quad b_1\in\Re(ACA-L_dAL_d),
\end{eqnarray}
where $b_1=(L_d-I)\beta$.\\
\noindent Proof:
The difference in MMSE is given by
\begin{eqnarray}
\Delta_1&=&MMSE(\hat\beta_{RMLE})-MMSE(\hat\beta_{RAULE})\cr
&=&ACA-L_dAL_d-b_1b_1'.
\end{eqnarray}
Since $ACA\geq 0$ and $L_dAL_d\geq 0$, making use of Lemma 2.4, under the given assumptions, one obtains $ACA-L_dAL_d\geq 0$. Then, the result follows by applying Lemma 2.2.

Now, necessary and sufficient
conditions for the new estimator to be superior to the RMLE in the MSE
sense, will be given.

\noindent \textbf{Theorem 3.4.} Under the assumptions of the logistic regression model (1), the RAULE is superior to the
RMLE in the MSE sense, if the biasing parameter $d$ satisfies the following inequality
 $$\frac{(\lambda_1+d)(\lambda_1+2-d)}{(1-d)^2}<\frac{max_i\alpha_i^2}{min_ia_{ii}},$$
where $\lambda_1\geq\ldots\geq\lambda_{p+1}$ are the ordered eigenvalues of $C$, $a_{ii}\geq 0$ is the ith diagonal element of the matrix $T'AT$ and $\alpha_i$ is the ith elements of $\beta'T$ for the orthogonal matrix $T$.\\
\noindent Proof. Consider the MSE difference
\begin{eqnarray}
\Delta_2&=&MSE(\hat\beta_{RMLE})-MSE(\hat\beta_{RAULE})\cr
&=&tr[MMSE(\hat\beta_{RMLE})-MMSE(\hat\beta_{RAULE})]\cr
&=&\sum_{i=1}^{p+1}\frac{1}{(\lambda_i+1)^4}[(\lambda_i+2-d)^2(\lambda_i+d)^2-(\lambda_i+1)^4]a_{ii}+(1-d)^4\alpha_i^2].
\end{eqnarray}
Differentiating $\Delta_2$ with respect to $d$, gives
\begin{eqnarray}
\frac{\partial\Delta_2}{\partial d}&=&\sum_{i=1}^{p+1}\frac{(1-d)^3}{(\lambda_i+1)^4}\left[\frac{4(\lambda_i+d)(\lambda_i+2-d)a_{ii}}{(1-d)^2}-4\alpha_i^2\right]\nonumber\\
&\leq&\sum_{i=1}^{p+1}\frac{4(1-d)^3}{(\lambda_i+1)^4}\left[\frac{(\lambda_1+d)(\lambda_1+2-d)min_ia_{ii}}{(1-d)^2}-max_i\alpha_i^2\right].
\end{eqnarray}
The result follows whenever
$$\frac{(\lambda_1+d)(\lambda_1+2-d)}{(1-d)^2}<\frac{max_i\alpha_i^2}{min_ia_{ii}}.$$

Now, we give comparison result between the RAULE and MLE in the MMSE sense.\\
\noindent\textbf{Theorem 3.5.} Under the assumptions of the logistic regression model (1), when
$\lambda_{max}(L_dAL_dC)\leq 1$,
RAULE is superior to the
 MLE in the MMSE sense if
and only if
\begin{eqnarray}
b_1'(C^{-1}-L_dAL_d)^+b_1\leq1,\quad b_1\in\Re(C^{-1}-L_dAL_d),
\end{eqnarray}
 where $b_1=(L_d-I)\beta$.\\
 \noindent Proof.
The difference in MMSE is given by
\begin{eqnarray}
\Delta_3&=&MMSE(\hat\beta_{MLE})-MMSE(\hat\beta_{RAULE})\cr
&=&C^{-1}-L_dAL_d-b_1b_1'.
\end{eqnarray}
Since $C^{-1}>0$ and $L_dAL_d\geq 0$, applying Lemma 2.3 yields $C^{-1}-L_dAL_d\geq 0$. Then the result follows applying  Lemma 2.2.

\noindent \textbf{Theorem 3.6.} Under the assumptions of the logistic regression model (1), when $$\frac{(\lambda_1+d)(\lambda_1+2-d)}{(1-d)^2}<\frac{max_i\alpha_i^2}{min_ia_{ii}},$$ the RAULE is superior to MLE in the sense of MSE criterion.

\noindent Proof. Consider the MSE difference given by
\begin{eqnarray}
\Delta_4&=&MSE(\hat\beta_{MLE})-MSE(\hat\beta_{RAULE})\cr
&=&\mbox{tr}[MMSE(\hat\beta_{MLE})-MMSE(\hat\beta_{RAULE})]\cr
&=&\mbox{tr}[MMSE(\hat\beta_{MLE})-MMSE(\hat\beta_{RMLE})\cr
&&+MMSE(\hat\beta_{RMLE})-MMSE(\hat\beta_{RAULE})].
\end{eqnarray}
Since
\begin{eqnarray*}
\mbox{tr}[MMSE(\hat\beta_{MLE})-MMSE(\hat\beta_{RMLE})]&=&\mbox{tr}[C^{-1}-ACA]\cr
&=&\mbox{tr}[C^{-1}-A]\cr
&=&\mbox{tr}[C^{-1}H'(HC^{-1}H')^{-1}HC^{-1}]\geq 0
\end{eqnarray*}
we have $\mbox{tr}[MMSE(\hat\beta_{RMLE})-MMSE(\hat\beta_{MAULE})]\geq0$. Hence $\Delta_4\geq 0$, using Theorem 3.4.

Finally we present the comparison between the RAULE and AULE in the MMSE sense.\\
\noindent\textbf{Theorem 3.7.} Under the assumptions of the logistic regression model (1), the RAULE is always superior to AULE in the MMSE sense.\\
 \noindent Proof.
Consider the following difference in MMSE
\begin{eqnarray}
\Delta_5&=&MMSE(\hat\beta_{AULE})-MMSE(\hat\beta_{RAULE})\cr
&=&L_dC^{-1}L_d-L_dAL_d\cr
&=&L_d(C^{-1}-A)L_d.
\end{eqnarray}
Making use of Lemma 2.1, $C^{-1}-A=C^{-1}H'(HC^{-1}H')^{-1}HC^{-1}\geq 0$ and the proof is complete.\\

\noindent\textbf{Corollary  3.1.} The RAULE is always superior to the
AULE in the MSE sense.

\section{Monte Carlo Simulation}
\subsection{The design of the simulation}
In this section, we conduct a Monte Carlo simulation to compare the performances of the MLE, AULE, RMLE and RAULE under different scenarios. Since we want to compare the estimators when there multicollinearity is present, the main factor of the simulation is the degree of correlation $\left(\rho^2\right)$ among the explanatory variables. Hence, following Gibbons (1981) and Kibria (2003), we use the following formula to generate correlated variables:
\begin{equation}
x_{ij}=(1-\rho^2)^{1/2}z_{ij}+\rho z_{ip},\ i=1,...,n, \
j=1,...,p,
\end{equation}
where $z_{ij}$s are independent standard normal pseudo-random numbers. We consider three different values of the degree of correlation corresponding to $0.9, 0.99$ and $0.999$.

The number of $n$ observations is generated using the Bernoulli distribution with parameter 
$\pi _i$ such that 
\begin{equation}
\pi_i=\frac{exp(x_i'\beta)}{1+exp(x_i'\beta)}, \ i=1,2,\ldots,n,
\end{equation}
for the dependent variable.

We fit logistic regression models having $p=4$ and $p=8$ explanatory variables. We consider the sample sizes $50, 100$ and $200$. The parameter values of $\beta$ are chosen so that $\beta'\beta=1$.

Following M{\aa}nsson et al. (2015), we use different restriction matrices, to capture the effect of imposing restrictions on estimators, as follows:\\

I)  $ H_1=\left(
\begin{array}{cccc}
1&0&-2&1\\
1&-1&1&-1\\
\end{array}
\right)$  with $h_1=(0,0)'$  when $p=4$ and

II) $H_2=\left(
\begin{array}{cccccccc}
1&0&-2&1&-3&1&1&1\\
1&1&0&1&-3&1&-2&1\\
\end{array}
\right)$  with $h_2=(0,0)'$  when $p=8$.\\
The simulation is replicated $2000$ times and estimated MSE is evaluated using
$$ \mbox{MSE}=\frac{\sum_{r=1}^{2000}(\tilde{\beta}_r-\beta)'(\tilde{\beta}_r-\beta)}{2000},     $$
where $\tilde{\beta}_r$  is any estimator considered in this study, in the rth repetition. 

\subsection{Results of the simulation}
We summarized the estimated MSE values of the estimators in Tables 1-6. It can be observed that an increase in the degree of correlation affects the MSE values of the estimators negatively. It is also observed that MSE has the worst performance (having the most MSE values) in all of the situations considered. The RMLE has always less MSE than that of MLE.

The performances of the AULE and RAULE depend on the parameter $d$. The smaller the value of $d$, the better performance. The RMLE has lesser MSE value than AULE when $d>0.4$ in almost all cases. 

Moreover, it can be deduced that our new estimator RAULE has always the least MSE value. Especially, the RAULE is the most robust option against the correlation when the value of $d$ is small.

From tables, it is realized that an increase in the sample size causes a decrease in the MSE values. Increasing the number of explanatory variables also affects the performances of the estimator negatively. Again the newly proposed RAULE is the most robust option for this situation.

\begin{table} \tiny
\center
   \caption{The estimated MSE values for different $d$ when $n=50$ and $p=4$}
{\begin{tabular}{@{}lcccccccccc}
\noalign{\smallskip}\hline
 & $d=0.1$ & $d=0.2$ &$d=0.3$& $d=0.4$& $d=0.5$ & $d=0.6$&$d=0.7$ &$d=0.8$ &$d=0.9$&$d=0.99$\\
\noalign{\smallskip}\hline
$\gamma=0.9$\\
MLE	&	14.6207	&	14.6207	&	14.6207	&	14.6207	&	14.6207	&	14.6207	&	14.6207	&	14.6207	&	14.6207	&	14.6207	\\
AULE	&	2.4696	&	4.0121	&	5.8178	&	7.7198	&	9.5737	&	11.2573	&	12.6704	&	13.7354	&	14.3966	&	14.6185	\\
RMLE	&	3.3376	&	3.3376	&	3.3376	&	3.3376	&	3.3376	&	3.3376	&	3.3376	&	3.3376	&	3.3376	&	3.3376	\\
RAULE	&	0.9280	&	1.2590	&	1.6281	&	2.0065	&	2.3691	&	2.6946	&	2.9658	&	3.1692	&	3.2950	&	3.3372	\\

$\gamma=0.99$\\
MLE	&	122.8774	&	122.8774	&	122.8774	&	122.8774	&	122.8774	&	122.8774	&	122.8774	&	122.8774	&	122.8774	&	122.8774	\\
AULE	&	6.3439	&	18.3859	&	34.4887	&	52.5983	&	70.9343	&	87.9906	&	102.5346	&	113.6079	&	120.5258	&	122.8537	\\
RMLE	&	17.6015	&	17.6015	&	17.6015	&	17.6015	&	17.6015	&	17.6015	&	17.6015	&	17.6015	&	17.6015	&	17.6015	\\
RAULE	&	1.4451	&	3.2122	&	5.4877	&	8.0031	&	10.5255	&	12.8580	&	14.8392	&	16.3440	&	17.2826	&	17.5983	\\

$\gamma=0.999$\\
MLE	&	2988.9881	&	2988.9881	&	2988.9881	&	2988.9881	&	2988.9881	&	2988.9881	&	2988.9881	&	2988.9881	&	2988.9881	&	2988.9881	\\
AULE	&	109.4131	&	389.4573	&	779.6334	&	1226.2873	&	1682.9183	&	2110.1799	&	2475.8795	&	2754.9782	&	2929.5914	&	2988.3912	\\
RMLE	&	279.7609	&	279.7609	&	279.7609	&	279.7609	&	279.7609	&	279.7609	&	279.7609	&	279.7609	&	279.7609	&	279.7609	\\
RAULE	&	10.7511	&	37.0637	&	73.5805	&	115.3140	&	157.9420	&	197.8068	&	231.9161	&	257.9424	&	274.2231	&	279.7052	\\
   \noalign{\smallskip}\hline
  \end{tabular}}
\end{table}

\begin{table}\tiny
   \caption{The estimated MSE values for different $d$ when $n=100$ and $p=4$}
{\begin{tabular}{@{}lcccccccccc}
\noalign{\smallskip}\hline
 & $d=0.1$ & $d=0.2$ &$d=0.3$& $d=0.4$& $d=0.5$ & $d=0.6$&$d=0.7$ &$d=0.8$ &$d=0.9$&$d=0.99$\\
\noalign{\smallskip}\hline
$\gamma=0.9$\\
MLE	&	5.7257	&	5.7257	&	5.7257	&	5.7257	&	5.7257	&	5.7257	&	5.7257	&	5.7257	&	5.7257	&	5.7257	\\
AULE	&	2.5062	&	3.0270	&	3.5549	&	4.0642	&	4.5327	&	4.9417	&	5.2757	&	5.5229	&	5.6746	&	5.7252	\\
RMLE	&	1.7583	&	1.7583	&	1.7583	&	1.7583	&	1.7583	&	1.7583	&	1.7583	&	1.7583	&	1.7583	&	1.7583	\\
RAULE	&	1.1074	&	1.2195	&	1.3292	&	1.4325	&	1.5259	&	1.6063	&	1.6714	&	1.7192	&	1.7484	&	1.7582	\\

$\gamma=0.99$\\
MLE	&	50.0384	&	50.0384	&	50.0384	&	50.0384	&	50.0384	&	50.0384	&	50.0384	&	50.0384	&	50.0384	&	50.0384	\\
AULE	&	4.2987	&	9.4622	&	15.9754	&	23.1044	&	30.2129	&	36.7629	&	42.3141	&	46.5238	&	49.1475	&	50.0294	\\
RMLE	&	9.8945	&	9.8945	&	9.8945	&	9.8945	&	9.8945	&	9.8945	&	9.8945	&	9.8945	&	9.8945	&	9.8945	\\
RAULE	&	1.5374	&	2.5562	&	3.7795	&	5.0856	&	6.3692	&	7.5410	&	8.5281	&	9.2737	&	9.7373	&	9.8929	\\

$\gamma=0.999$\\
MLE	&	461.9003	&	461.9003	&	461.9003	&	461.9003	&	461.9003	&	461.9003	&	461.9003	&	461.9003	&	461.9003	&	461.9003	\\
AULE	&	18.9701	&	62.9223	&	123.3250	&	192.0703	&	262.1311	&	327.5616	&	383.4969	&	426.1535	&	452.8284	&	461.8091	\\
RMLE	&	124.0645	&	124.0645	&	124.0645	&	124.0645	&	124.0645	&	124.0645	&	124.0645	&	124.0645	&	124.0645	&	124.0645	\\
RAULE	&	5.5594	&	17.3608	&	33.5400	&	51.9347	&	70.6708	&	88.1627	&	103.1129	&	114.5124	&	121.6405	&	124.0402	\\

   \noalign{\smallskip}\hline
  \end{tabular}}
\end{table}

\begin{table}\tiny
   \caption{The estimated MSE values for different $d$ when $n=200$ and $p=4$}
{\begin{tabular}{@{}lcccccccccc}
\noalign{\smallskip}\hline
 & $d=0.1$ & $d=0.2$ &$d=0.3$& $d=0.4$& $d=0.5$ & $d=0.6$&$d=0.7$ &$d=0.8$ &$d=0.9$&$d=0.99$\\
\noalign{\smallskip}\hline
$\gamma=0.9$\\
MLE	&	2.4237	&	2.4237	&	2.4237	&	2.4237	&	2.4237	&	2.4237	&	2.4237	&	2.4237	&	2.4237	&	2.4237	\\
AULE	&	1.8950	&	1.9963	&	2.0899	&	2.1743	&	2.2481	&	2.3100	&	2.3592	&	2.3949	&	2.4165	&	2.4237	\\
RMLE	&	1.1825	&	1.1825	&	1.1825	&	1.1825	&	1.1825	&	1.1825	&	1.1825	&	1.1825	&	1.1825	&	1.1825	\\
RAULE	&	1.0220	&	1.0532	&	1.0818	&	1.1074	&	1.1297	&	1.1484	&	1.1632	&	1.1739	&	1.1803	&	1.1825	\\

$\gamma=0.99$\\
MLE	&	18.7164	&	18.7164	&	18.7164	&	18.7164	&	18.7164	&	18.7164	&	18.7164	&	18.7164	&	18.7164	&	18.7164	\\
AULE	&	3.4884	&	5.5229	&	7.8305	&	10.2190	&	12.5219	&	14.5982	&	16.3327	&	17.6356	&	18.4430	&	18.7137	\\
RMLE	&	4.8649	&	4.8649	&	4.8649	&	4.8649	&	4.8649	&	4.8649	&	4.8649	&	4.8649	&	4.8649	&	4.8649	\\
RAULE	&	1.4912	&	1.9655	&	2.4871	&	3.0174	&	3.5228	&	3.9750	&	4.3508	&	4.6321	&	4.8061	&	4.8643	\\

$\gamma=0.999$\\
MLE	&	168.3497	&	168.3497	&	168.3497	&	168.3497	&	168.3497	&	168.3497	&	168.3497	&	168.3497	&	168.3497	&	168.3497	\\
AULE	&	8.5172	&	24.9718	&	47.0305	&	71.8658	&	97.0272	&	120.4414	&	140.4116	&	155.6186	&	165.1198	&	168.3173	\\
RMLE	&	52.5252	&	52.5252	&	52.5252	&	52.5252	&	52.5252	&	52.5252	&	52.5252	&	52.5252	&	52.5252	&	52.5252	\\
RAULE	&	3.0800	&	8.1803	&	15.0086	&	22.6921	&	30.4739	&	37.7139	&	43.8883	&	48.5895	&	51.5267	&	52.5151	\\

   \noalign{\smallskip}\hline
  \end{tabular}}
\end{table}

\begin{table}\tiny
   \caption{The estimated MSE values for different $d$ when $n=50$ and $p=8$}
{\begin{tabular}{@{}lcccccccccc}
\noalign{\smallskip}\hline
 & $d=0.1$ & $d=0.2$ &$d=0.3$& $d=0.4$& $d=0.5$ & $d=0.6$&$d=0.7$ &$d=0.8$ &$d=0.9$&$d=0.99$\\
\noalign{\smallskip}\hline
$\gamma=0.9$\\
MLE	&	141.9664	&	141.9664	&	141.9664	&	141.9664	&	141.9664	&	141.9664	&	141.9664	&	141.9664	&	141.9664	&	141.9664	\\
AULE	&	8.3150	&	22.3253	&	40.8815	&	61.6607	&	82.6499	&	102.1456	&	118.7542	&	131.3918	&	139.2840	&	141.9395	\\
RMLE	&	42.1899	&	42.1899	&	42.1899	&	42.1899	&	42.1899	&	42.1899	&	42.1899	&	42.1899	&	42.1899	&	42.1899	\\
RAULE	&	3.5096	&	7.7610	&	13.2181	&	19.2412	&	25.2760	&	30.8532	&	35.5891	&	39.1851	&	41.4281	&	42.1823	\\

$\gamma=0.99$\\
MLE	&	1270.2249	&	1270.2249	&	1270.2249	&	1270.2249	&	1270.2249	&	1270.2249	&	1270.2249	&	1270.2249	&	1270.2249	&	1270.2249	\\
AULE	&	48.8557	&	168.8119	&	334.8236	&	524.3276	&	717.7695	&	898.6031	&	1053.2906	&	1171.3026	&	1245.1183	&	1269.9726	\\
RMLE	&	368.8009	&	368.8009	&	368.8009	&	368.8009	&	368.8009	&	368.8009	&	368.8009	&	368.8009	&	368.8009	&	368.8009	\\
RAULE	&	15.1951	&	50.3383	&	98.5839	&	153.4681	&	209.3890	&	261.6062	&	306.2416	&	340.2786	&	361.5626	&	368.7282	\\

$\gamma=0.999$\\
MLE	&	46865.2930	&	46865.2930	&	46865.2930	&	46865.2930	&	46865.2930	&	46865.2930	&	46865.2930	&	46865.2930	&	46865.2930	&	46865.2930	\\
AULE	&	1693.7814	&	6076.5807	&	12192.7161	&	19198.8276	&	26364.0035	&	33069.7799	&	38810.1409	&	43191.5188	&	45932.7934	&	46855.9216	\\
RMLE	&	12492.1038	&	12492.1038	&	12492.1038	&	12492.1038	&	12492.1038	&	12492.1038	&	12492.1038	&	12492.1038	&	12492.1038	&	12492.1038	\\
RAULE	&	452.1749	&	1620.7160	&	3251.0570	&	5118.4698	&	7028.1898	&	8815.4162	&	10345.3121	&	11513.0044	&	12243.5835	&	12489.6063	\\

   \noalign{\smallskip}\hline
  \end{tabular}}
\end{table}
\begin{table}\tiny
   \caption{The estimated MSE values for different $d$ when $n=100$ and $p=8$}
{\begin{tabular}{@{}lcccccccccc}
\noalign{\smallskip}\hline
 & $d=0.1$ & $d=0.2$ &$d=0.3$& $d=0.4$& $d=0.5$ & $d=0.6$&$d=0.7$ &$d=0.8$ &$d=0.9$&$d=0.99$\\
\noalign{\smallskip}\hline
$\gamma=0.9$\\
MLE	&	26.2009	&	26.2009	&	26.2009	&	26.2009	&	26.2009	&	26.2009	&	26.2009	&	26.2009	&	26.2009	&	26.2009	\\
AULE	&	4.6742	&	7.4254	&	10.6324	&	14.0028	&	17.2834	&	20.2599	&	22.7567	&	24.6377	&	25.8052	&	26.1969	\\
RMLE	&	11.8713	&	11.8713	&	11.8713	&	11.8713	&	11.8713	&	11.8713	&	11.8713	&	11.8713	&	11.8713	&	11.8713	\\
RAULE	&	2.7895	&	4.0017	&	5.3773	&	6.8017	&	8.1754	&	9.4141	&	10.4490	&	11.2264	&	11.7082	&	11.8697	\\

$\gamma=0.99$\\
MLE	&	309.5121	&	309.5121	&	309.5121	&	309.5121	&	309.5121	&	309.5121	&	309.5121	&	309.5121	&	309.5121	&	309.5121	\\
AULE	&	14.7981	&	44.9405	&	85.5269	&	131.3111	&	177.7461	&	220.9848	&	257.8793	&	285.9813	&	303.5420	&	309.4521	\\
RMLE	&	141.9314	&	141.9314	&	141.9314	&	141.9314	&	141.9314	&	141.9314	&	141.9314	&	141.9314	&	141.9314	&	141.9314	\\
RAULE	&	7.5535	&	21.5523	&	40.1706	&	61.0585	&	82.1795	&	101.8102	&	118.5408	&	131.2744	&	139.2280	&	141.9042	\\

$\gamma=0.999$\\
MLE	&	3051.3774	&	3051.3774	&	3051.3774	&	3051.3774	&	3051.3774	&	3051.3774	&	3051.3774	&	3051.3774	&	3051.3774	&	3051.3774	\\
AULE	&	113.5289	&	400.3017	&	798.8422	&	1254.5856	&	1720.2427	&	2155.7995	&	2528.5172	&	2812.9323	&	2990.8567	&	3050.7692	\\
RMLE	&	1382.3472	&	1382.3472	&	1382.3472	&	1382.3472	&	1382.3472	&	1382.3472	&	1382.3472	&	1382.3472	&	1382.3472	&	1382.3472	\\
RAULE	&	52.1509	&	182.3238	&	362.9193	&	569.2861	&	780.0596	&	977.1620	&	1145.8025	&	1274.4772	&	1354.9689	&	1382.0721	\\

   \noalign{\smallskip}\hline
  \end{tabular}}
\end{table}
\begin{table}\tiny
   \caption{The estimated MSE values for different $d$ when $n=200$ and $p=8$}
{\begin{tabular}{@{}lcccccccccc}
\noalign{\smallskip}\hline
 & $d=0.1$ & $d=0.2$ &$d=0.3$& $d=0.4$& $d=0.5$ & $d=0.6$&$d=0.7$ &$d=0.8$ &$d=0.9$&$d=0.99$\\
\noalign{\smallskip}\hline
$\gamma=0.9$\\
MLE	&	6.6427	&	6.6427	&	6.6427	&	6.6427	&	6.6427	&	6.6427	&	6.6427	&	6.6427	&	6.6427	&	6.6427	\\
AULE	&	3.9466	&	4.4224	&	4.8819	&	5.3103	&	5.6947	&	6.0240	&	6.2895	&	6.4841	&	6.6028	&	6.6423	\\
RMLE	&	4.1607	&	4.1607	&	4.1607	&	4.1607	&	4.1607	&	4.1607	&	4.1607	&	4.1607	&	4.1607	&	4.1607	\\
RAULE	&	2.6886	&	2.9525	&	3.2052	&	3.4393	&	3.6483	&	3.8268	&	3.9703	&	4.0752	&	4.1392	&	4.1604	\\

$\gamma=0.99$\\
MLE	&	74.7216	&	74.7216	&	74.7216	&	74.7216	&	74.7216	&	74.7216	&	74.7216	&	74.7216	&	74.7216	&	74.7216	\\
AULE	&	7.8025	&	15.6628	&	25.3268	&	35.7714	&	46.1097	&	55.5915	&	63.6028	&	69.6663	&	73.4407	&	74.7088	\\
RMLE	&	39.4314	&	39.4314	&	39.4314	&	39.4314	&	39.4314	&	39.4314	&	39.4314	&	39.4314	&	39.4314	&	39.4314	\\
RAULE	&	5.0776	&	9.2535	&	14.2768	&	19.6453	&	24.9242	&	29.7452	&	33.8072	&	36.8760	&	38.7842	&	39.4249	\\

$\gamma=0.999$\\
MLE	&	815.6036	&	815.6036	&	815.6036	&	815.6036	&	815.6036	&	815.6036	&	815.6036	&	815.6036	&	815.6036	&	815.6036	\\
AULE	&	33.8851	&	111.8311	&	218.6004	&	339.9454	&	463.5184	&	578.8710	&	677.4549	&	752.6211	&	799.6205	&	815.4430	\\
RMLE	&	411.3009	&	411.3009	&	411.3009	&	411.3009	&	411.3009	&	411.3009	&	411.3009	&	411.3009	&	411.3009	&	411.3009	\\
RAULE	&	18.0382	&	57.6046	&	111.4736	&	172.5361	&	234.6310	&	292.5449	&	342.0125	&	379.7161	&	403.2862	&	411.2203	\\

   \noalign{\smallskip}\hline
  \end{tabular}}
\end{table}
\section{Application}
In this section, we present a real data application in order to show the benefit of using the newly proposed RAULE. For our purpose, we used a data set available at an official web page of Statistics Sweden (http://www.scb.se/). The observations are 83 municipalities which are the urban regions belonging to the functional analysis regions Stockholm, Malmo and Goteborg. Asar and Genc (2015) and Mansson et al. (2012) also analyzed similar data sets. We model the data using a binary logistic regression model such that the dependent variable is coded as 1 if there is an increase in the pupation and 0 if there is a decrease. The dependent variable is explained by the following explanatory variables:

X1: The population,

X2: The number of unemployed people,

X3: The number of newly constructed buildings,

X4: The number of bankrupt firms.\\
The correlation matrix of this data set is given in Table 7. It is
observed from Table 7 that all the bivariate correlations between
the explanatory variables are larger than 0.95. The condition
number, being a measure of the degree of multicollinearity, is
computed according to  $\kappa=\sqrt{\lambda_{max}/\lambda_{min}}$
which shows that there exists a severe multicollinearity problem.

\begin{table}\scriptsize
   \caption{The correlation matrix of the design matrix}
   \begin{center}
{\begin{tabular}{@{}lccccc}
\noalign{\smallskip}\hline
 & X1 & X2 &X3& X4\\
\noalign{\smallskip}\hline
  X1  &   1.0000 &    0.9937  & 0.9707 &    0.9514\\
  X2  &   0.9937 &    1.0000  & 0.9527 &    0.9222\\
  X3  &   0.9707 &    0.9527  & 1.0000 &    0.9765\\
  X4  &   0.9514 &    0.9222  & 0.9765 &    1.0000\\

   \noalign{\smallskip}\hline
  \end{tabular}}
  \end{center}
\end{table}

\begin{table}\scriptsize
   \caption{The estimated theoretical MSE value of AULE}
   \begin{center}
{\begin{tabular}{@{}lccccccc}
\noalign{\smallskip}\hline
 d&$\beta_1$ & $\beta_2$  &$\beta_3$ &$\beta_4$ &MSE\\
\noalign{\smallskip}\hline
0.1&  5.0078&   -3.1222&    0.8875& -2.0282 &70.9595\\
0.2&9.2699& -6.1203&    1.5128  &-3.9040&   248.9311\\
0.3&13.0305&    -8.7656&    2.0645  &-5.5592&   496.2849\\
0.4&16.2896 &-11.0583&  2.5427  &-6.9937&   779.1508\\
0.5&19.0474&    -12.9982&   2.9473& -8.2075 &1068.1749\\
0.6&21.3038&    -14.5854&   3.2784& -9.2006&    1338.5192\\
0.7&23.0587&    -15.8199&   3.5359  &-9.9730&   1569.8617\\
0.8&24.3122&    -16.7017&   3.7198& -10.5248&   1746.3963\\
0.9&25.0644&    -17.2307&   3.8301& -10.8558&   1856.8330\\
0.99&25.3126&   -17.4053&   3.8666& -10.9650&   1894.0204\\

   \noalign{\smallskip}\hline
  \end{tabular}}
  \end{center}
\end{table}

\begin{table}\scriptsize
   \caption{The estimated theoretical MSE value of RAULE}
   \begin{center}
{\begin{tabular}{@{}lccccccc}
\noalign{\smallskip}\hline
 d&$\beta_1$ & $\beta_2$  &$\beta_3$ &$\beta_4$ &MSE\\
\noalign{\smallskip}\hline
0.1&0.5427  &0.5772 &0.4917 &-0.8105&   4.3052\\
0.2&0.8249  &0.8521 &0.7846 &-1.5974    &12.2573\\
0.3&1.0739& 1.0947  &1.0430&    -2.2917&    23.0005\\
0.4&1.2896& 1.3050  &1.2670&    -2.8934 &35.1374\\
0.5&1.4722& 1.4829& 1.4565& -3.4026&    47.4568\\
0.6&1.6216& 1.6284  &1.6115 &-3.8192&   58.9338\\
0.7&1.7378& 1.7416& 1.7321& -4.1432&    68.7299\\
0.8&1.8208& 1.8225& 1.8183& -4.3747&    76.1930\\
0.9&1.8706& 1.8710& 1.8699& -4.5135&    80.8572\\
0.99&1.8870&    1.8870& 1.8870& -4.5594&    82.4270\\

   \noalign{\smallskip}\hline
  \end{tabular}}
  \end{center}
\end{table}

\begin{figure}
\begin{center}{
{\includegraphics[width=10cm,height=7cm]{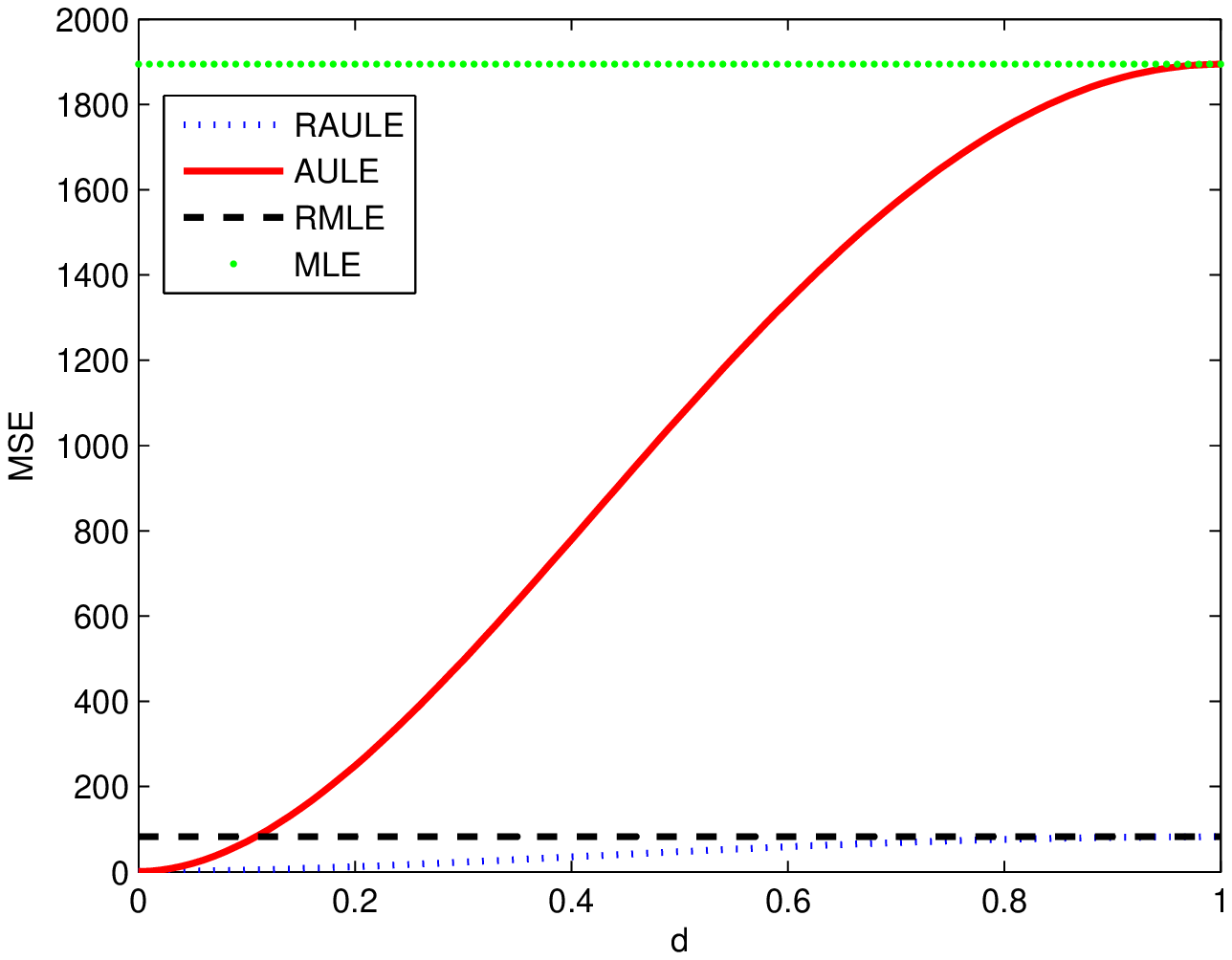}}}
 \caption{ The estimated MSE values of the estimators when $0\leq d\leq 1$}
\end{center}
\end{figure}
We computed the estimated theoretical MSE values along with the coefficients of the estimators. The MSE of MLE is 1894.398 which is the largest among others. The coefficients of MLE are 25.3151, -17.4071, 3.8669 and -10.9661. The MSE of RMLE is 82.4430 and the coefficients of RMLE are 1.8872, 1.8872, 1.8872 and -4.5598. MSE values and the coefficients of AULE and RAULE are given in Tables 8 and 9 respectively for different biasing parameter $d$ varying from zero to one.

According to Tables 8-9, the RAULE has the least MSE value for all
values of the parameter $d$. When $d=1$, the MSE of RAULE becomes equal
to that of RMLE as expected. The MSE of AULE is lower than that of RMLE
when $d<0.12$. We also provided the plot of the MSE versus $d$ in
Figure 1. According to Figure 1, the RAULE has the best performance
among other estimators.

\section{Conclusions}
\noindent In this paper, we proposed a restricted almost unbiased Liu estimator in the logistic regression model. Then, its performance compared with other competitors including the MLE, AULE, RMLE in the logistic regression model through providing some theorems, a Monte Carlo simulation as well as a real example. We concluded that our proposed estimator is superior compared to all others in the sense of having smaller MSE value.




\bibliographystyle{elsarticle-num}

\end{document}